\documentclass[12pt,leqno,draft]{article}
\usepackage{amsfonts}
\usepackage{amsmath, dsfont, amsthm, amsfonts, amssymb, color}
\usepackage{mathrsfs}
\usepackage{color}
\usepackage{stmaryrd}
\newtheorem{thm}{Theorem}[section]

\newtheorem{lem}[thm]{Lemma}

\newtheorem{rem}[thm]{Remark}
\theoremstyle{definition}

\newcommand{\scr}[1]{\mathscr #1}
\definecolor{wco}{rgb}{0.5,0.2,0.3}

\numberwithin{equation}{section} \theoremstyle{remark}

\newcommand{\ua}{\uparrow}

\title{{\bf    Probability  Distance Estimates  Between   Diffusion Processes and Applications to Singular McKean-Vlasov SDEs}\footnote{ Supported in
 part by   the National Key R\&D Program of China (2022YFA1006000, 2020YFA0712900) and   NNSFC (12271398, 11921001).}
 }
\author{
{\bf   Xing Huang $^{a)}$,  Panpan Ren $^{b)}$, Feng-Yu Wang $^{a)}$  }\\
\footnotesize{ a)Center for Applied Mathematics, Tianjin
University, Tianjin 300072, China}\\
 \footnotesize{ b)Department of Mathematics,
City University of Hong Kong, Tat Chee Avenue, Hong Kong,  China}\\
\footnotesize{  xinghuang@tju.edu.cn, panparen@cityu.edu.hk, wangfy@tju.edu.cn}}
\begin{document}
\allowdisplaybreaks
\def\R{\mathbb R}  \def\ff{\frac} \def\ss{\sqrt} \def\B{\mathbf
B}
\def\N{\mathbb N} \def\kk{\kappa} \def\m{{\bf m}}
\def\ee{\varepsilon}\def\ddd{D^*}
\def\dd{\delta} \def\DD{\Delta} \def\vv{\varepsilon} \def\rr{\rho}
\def\<{\langle} \def\>{\rangle} \def\GG{\Gamma} \def\gg{\gamma}
  \def\nn{\nabla} \def\pp{\partial} \def\E{\mathbb E}
\def\d{\text{\rm{d}}} \def\bb{\beta} \def\aa{\alpha} \def\D{\scr D}
  \def\si{\sigma} \def\ess{\text{\rm{ess}}}
\def\beg{\begin} \def\beq{\begin{equation}}  \def\F{\scr F}
\def\Ric{\text{\rm{Ric}}} \def\Hess{\text{\rm{Hess}}}
\def\e{\text{\rm{e}}} \def\ua{\underline a} \def\OO{\Omega}  \def\oo{\omega}
 \def\tt{\tilde} \def\Ric{\text{\rm{Ric}}}
\def\cut{\text{\rm{cut}}} \def\P{\mathbb P} \def\ifn{I_n(f^{\bigotimes n})}
\def\C{\scr C}   \def\G{\scr G}   \def\aaa{\mathbf{r}}     \def\r{r}
\def\gap{\text{\rm{gap}}} \def\prr{\pi_{{\bf m},\varrho}}  \def\r{\mathbf r}
\def\Z{\mathbb Z} \def\vrr{\varrho} \def\ll{\lambda}
\def\L{\scr L}\def\Tt{\tt} \def\TT{\tt}\def\II{\mathbb I}
\def\i{{\rm in}}\def\Sect{{\rm Sect}}  \def\H{\mathbb H}
\def\M{\scr M}\def\Q{\mathbb Q} \def\texto{\text{o}} \def\LL{\Lambda}
\def\Rank{{\rm Rank}} \def\B{\scr B} \def\i{{\rm i}} \def\HR{\hat{\R}^d}
\def\to{\rightarrow}\def\l{\ell}\def\iint{\int}
\def\EE{\scr E}\def\no{\nonumber}
\def\A{\scr A}\def\V{\mathbb V}\def\osc{{\rm osc}}
\def\BB{\scr B}\def\Ent{{\rm Ent}}\def\3{\triangle}\def\H{\scr H}
\def\U{\scr U}\def\8{\infty}\def\1{\lesssim}\def\HH{\mathrm{H}}
 \def\T{\scr T}
 \def\R{\mathbb R}  \def\ff{\frac} \def\ss{\sqrt} \def\B{\mathbf
B} \def\W{\mathbb W}
\def\N{\mathbb N} \def\kk{\kappa} \def\m{{\bf m}}
\def\ee{\varepsilon}\def\ddd{D^*}
\def\dd{\delta} \def\DD{\Delta} \def\vv{\varepsilon} \def\rr{\rho}
\def\<{\langle} \def\>{\rangle} \def\GG{\Gamma} \def\gg{\gamma}
  \def\nn{\nabla} \def\pp{\partial} \def\E{\mathbb E}
\def\d{\text{\rm{d}}} \def\bb{\beta} \def\aa{\alpha} \def\D{\scr D}
  \def\si{\sigma} \def\ess{\text{\rm{ess}}}
\def\beg{\begin} \def\beq{\begin{equation}}  \def\F{\scr F}
\def\Ric{\text{\rm{Ric}}} \def\Hess{\text{\rm{Hess}}}
\def\e{\text{\rm{e}}} \def\ua{\underline a} \def\OO{\Omega}  \def\oo{\omega}
 \def\tt{\tilde} \def\Ric{\text{\rm{Ric}}}
\def\cut{\text{\rm{cut}}} \def\P{\mathbb P} \def\ifn{I_n(f^{\bigotimes n})}
\def\C{\scr C}      \def\aaa{\mathbf{r}}     \def\r{r}
\def\gap{\text{\rm{gap}}} \def\prr{\pi_{{\bf m},\varrho}}  \def\r{\mathbf r}
\def\Z{\mathbb Z} \def\vrr{\varrho} \def\ll{\lambda}
\def\L{\scr L}\def\Tt{\tt} \def\TT{\tt}\def\II{\mathbb I}
\def\i{{\rm in}}\def\Sect{{\rm Sect}}  \def\H{\mathbb H}
\def\M{\scr M}\def\Q{\mathbb Q} \def\texto{\text{o}} \def\LL{\Lambda}
\def\Rank{{\rm Rank}} \def\B{\scr B} \def\i{{\rm i}} \def\HR{\hat{\R}^d}
\def\to{\rightarrow}\def\l{\ell}\def\iint{\int}
\def\EE{\scr E}\def\Cut{{\rm Cut}}
\def\A{\scr A} \def\Lip{{\rm Lip}}
\def\BB{\scr B}\def\Ent{{\rm Ent}}\def\L{\scr L}
\def\R{\mathbb R}  \def\ff{\frac} \def\ss{\sqrt} \def\B{\mathbf
B}
\def\N{\mathbb N} \def\kk{\kappa} \def\m{{\bf m}}
\def\dd{\delta} \def\DD{\Delta} \def\vv{\varepsilon} \def\rr{\rho}
\def\<{\langle} \def\>{\rangle} \def\GG{\Gamma} \def\gg{\gamma}
  \def\nn{\nabla} \def\pp{\partial} \def\E{\mathbb E}
\def\d{\text{\rm{d}}} \def\bb{\beta} \def\aa{\alpha} \def\D{\scr D}
  \def\si{\sigma} \def\ess{\text{\rm{ess}}}
\def\beg{\begin} \def\beq{\begin{equation}}  \def\F{\scr F}
\def\Ric{\text{\rm{Ric}}} \def\Hess{\text{\rm{Hess}}}
\def\e{\text{\rm{e}}} \def\ua{\underline a} \def\OO{\Omega}  \def\oo{\omega}
 \def\tt{\tilde} \def\Ric{\text{\rm{Ric}}}
\def\cut{\text{\rm{cut}}} \def\P{\mathbb P} \def\ifn{I_n(f^{\bigotimes n})}
\def\C{\scr C}      \def\aaa{\mathbf{r}}     \def\r{r}
\def\gap{\text{\rm{gap}}} \def\prr{\pi_{{\bf m},\varrho}}  \def\r{\mathbf r}
\def\Z{\mathbb Z} \def\vrr{\varrho} \def\ll{\lambda}
\def\L{\scr L}\def\Tt{\tt} \def\TT{\tt}\def\II{\mathbb I}
\def\i{{\rm in}}\def\Sect{{\rm Sect}}  \def\H{\mathbb H}
\def\M{\scr M}\def\Q{\mathbb Q} \def\texto{\text{o}} \def\LL{\Lambda}
\def\Rank{{\rm Rank}} \def\B{\scr B} \def\i{{\rm i}} \def\HR{\hat{\R}^d}
\def\to{\rightarrow}\def\l{\ell}
\def\8{\infty}\def\I{1}\def\U{\scr U}\def\beq{\begin{equation}}
\maketitle

\begin{abstract} The $L^k$-Wasserstein distance $\W_k (k\ge 1)$ and the probability distance $\W_\psi$ induced by a concave function $\psi$,  are estimated      between different  diffusion processes with singular coefficients. As applications, the well-posedness, probability distance estimates and the log-Harnack inequality are derived for   McKean-Vlasov SDEs with multiplicative
distribution dependent noise, where the coefficients are singular in time-space variables and  $(\W_k+\W_\psi)$-Lipschitz continuous in the  distribution variable.
This improves existing results derived in the literature under  the $\W_k$-Lipschitz or derivative conditions in the  distribution variable.
 \end{abstract} \noindent
\vskip 2cm
\section{Introduction}

 Let $T>0$, and let $\Xi$ be the space of $(a,b)$, where
 $$b: [0,T]\times\R^d\to\R^d,\ \ \ a: [0,T]\times \R^d\to \R^d\otimes\R^d$$
 are measurable, and for any $(t,x)\in [0,T]\times \R^d,$ $a(t,x)$ is positive definite. For any $(a,b)\in \Xi,$ consider the  time dependent second order differential operator  on $\R^d$:
 $$L_{t}^{a,b}:=  {\rm tr}\{a(t,\cdot) \nn^2\}+ b(t,\cdot)\cdot\nn,\ \ t\in [0,T].$$
Let $(a_i,b_i)\in \Xi, i=1,2,$ such that for any $s\in [0,T)$, each $(L_t^{a_i,b_i})_{t\in [s,T]}$ generates a unique diffusion process $(X_{s,t}^{i,x})_{(t,x)\in [s,T]\times\R^d}$ on $\R^d$  with $X_{s,s}^{i,x}=x.$  Let
$$P_{s,t}^{i,x}:=\L_{X_{s,t}^{i,x}}$$ be the distribution of $X_{s,t}^{i,x}.$
When $s=0$, we simply denote
 $$X_{0,t}^{i,x}=X_t^{i,x},\ \ \ P_{0,t}^{i,x}= P_t^{i,x}.$$
If the initial value is random with distributions $\gg\in \scr P,$ where $\scr P$ is the set of all probability measures on $\R^d$, we denote the diffusion process  by $X_{s,t}^{i,\gamma}$, which has distribution
\beq\label{MK} P_{s,t}^{i,\gg} =\int_{\R^d} P_{s,t}^{i,x} \gg (\d x),\ \ \ i=1,2, \ 0\leq s\leq t\leq T.\end{equation}

By developing the bi-coupling argument and using an entropy inequality due to \cite{BRS},    the relative entropy
 $$\Ent(P_{s,t}^{1,\gg} |P_{s,t}^{2,\tt\gg}):=\int_{\R^d} \Big(\log \ff{\d P_{s,t}^{1,\gg} }{\d P_{s,t}^{2,\tt\gg}}\Big)\d P_{s,t}^{1,\gg},\ \ 0\leq s<t\leq T, \gamma,\tilde{\gamma}\in\scr P$$ is estimated in \cite{23RW}, and as an application, the log-Haranck inequality is established for McKean-Vlasov SDEs with multiplicative distribution dependent noise, where the drift is Dini continuous in the spatial variable $x$,  and the diffusion coefficient   is Lipschitz continuous  in $x$ and the distribution variable  with respect to
 $\W_2$.

 In this paper, we estimate a weighted variational distance   between $P_t^{1,\gg}$ and $P_t^{2,\tt\gg}$ for diffusion processes with singular coefficients,  and apply to the study of singular McKean-Vlasov SDEs with multiplicative distribution dependent noise, so that existing results in the literature are considerably extended.

Consider the class
$$ \scr A:= \big\{\psi: [0,\infty)\to  [0,\infty) \  \text{is\ increasing\ and\ concave,\ } \psi(r)>0\ \text{for}\ r>0\big\}. $$
For any $\psi\in \scr A,$ the $\psi$-continuity modulus of  a function $f$ on $\R^d$ is
$$[f]_\psi:=\sup_{x\ne y} \ff{|f(x)-f(y)|}{\psi(|x-y|)}.$$
Then
$$\scr P_\psi:=\bigg\{\mu\in \scr P:\
 \|\mu\|_\psi:=\int_{\R^d}\psi(|x|)\mu(\d x)<\infty\bigg\}$$ is a complete metric space under the
 distance $\W_\psi$ induced by $\psi$:
$$\W_\psi(\mu,\nu):=\sup_{[f]_\psi\le 1} |\mu(f)-\nu(f)|,$$
where $\mu(f):=\int_{\R^d}f\d\mu$ for $f\in L^1(\mu).$
In particular,   $\W_\psi=\W_1$ is the $L^1$-Wasserstein distance if $\psi(r)=r$, while $\W_\psi$ with  $\psi\equiv 2$  reduces to the total variational distance
$$\|\mu-\nu\|_{var} :=\sup_{|f|\le 1} |\mu(f)-\nu(f)|.$$
For any $k>0,$  the $L^k$-Wasserstein distance is
$$\W_k(\mu,\nu):=\inf_{\pi\in \C(\mu,\nu)} \bigg(\int_{\R^d\times\R^d} |x-y|^k\pi(\d x,\d y)\bigg)^{\ff 1 {1\lor k}},$$
where $\C(\mu,\nu)$ is the set of couplings for $\mu$ and $\nu$.  Then
$$\scr P_k:=\big\{\mu\in \scr P:\ \mu(|\cdot|^k)<\infty\big\}$$ is a Polish space under $\W_k.$ Since $\psi$ has at most linear growth, we have $ \scr P_k\subset \scr P_\psi, $
and $\scr P_k$ is complete under $\W_\psi+\W_k.$

To characterize  the singularity of coefficients in time-space variables, we recall some functional spaces introduced in \cite{XXZZ}.  For any $p\ge 1$, $L^p(\R^d)$ is the class of   measurable  functions $f$ on $\R^d$ such that
 $$\|f\|_{L^p(\R^d)}:=\bigg(\int_{\R^d}|f(x)|^p\d x\bigg)^{\ff 1 p}<\infty.$$
For any $p,q>1$ and a measurable function  $f$ on $[0,T]\times\R^d$, let
$$\|f\|_{\tt L_q^p(s,t)}:= \sup_{z\in \R^d}\bigg( \int_{s}^{t} \|1_{B(z,1)}f_r\|_{L^p(\R^d)}^q\d r\bigg)^{\ff 1 q},$$
where $B(z,1):= \{x\in\R^d: |x-z|\le 1\}$. When $s=0,$ we simply   denote $\|\cdot\|_{\tt L_q^p(t)}= \|\cdot\|_{\tt L_q^p(0,t)}$.
Let
$$\scr K:=\Big\{(p,q): p,q\in  (2,\infty),\  \ff d p+\ff 2 q<1\Big\}.$$  Let  $\|\cdot\|_\infty$ be the uniform norm, and for any function $f$ on $[0,T]\times \R^d,$ let
$$\|f\|_{t,\infty}:=\sup_{x\in\R^d}|f(t,x)|,\ \ \|f\|_{r\to t,\infty}:=\sup_{s\in [r,t]}\|f\|_{s,\infty},\ \ \ 0\le r\le t\le T.$$
  We make the following assumptions for the coefficients $(a,b)\in\Xi$, where $\nn$ is the gradient operator on $\R^d$.

\begin{enumerate}
\item[$(A^{a,b})$] There exist    constants $\aa\in (0,1], K>1, l\in \mathbb N$ and   $\{(p_i,q_i)\}_{0\le i\le l}\subset \scr K$ such that the following conditions hold.
\item[$(1)$]
 $\|a\|_\infty\lor \|a^{-1}\|_\infty \leq K,$ and
\beq\label{AA}    \|a (t, x) -a(t,y)\|\le K |x-y|^\aa,\ \ \ t\in [0,T], x,y\in \R^d.\end{equation}
 Moreover, there exist $\{1\le f_i\}_{1\le i\le l} $ with $\sum_{i=1}^l \|f_i\|_{\tt L_{q_i}^{p_i}(T)} \le K,$ such that
 $$\|\nn a\|\le \sum_{i=1}^l f_i.$$
 \item[$(2)$] $b$ has a decomposition $b=b^{(0)}+b^{(1)}$ such that
 $$\sup_{t\in [0,T]}|b^{(1)}(t,0)|+ \|\nn b^{(1)}\|_\infty + \|b^{(0)}\|_{\tt L_{q_0}^{p_0}(T)}\le K.$$  \end{enumerate}

 Let $\si(t,x):=\ss{2 a(t,x)}$, and let $W_t$ be a $d$-dimensional Brownian motion on a probability basis $(\OO,\F, \{\F_t\}_{t\in [0,T]},\P)$. By \cite[Theorem 2.1]{Ren} for $V(x):=1+|x|^2$, see also \cite{XXZZ} or \cite{ZY},
 under $(A^{a,b})$, for any $(s,x)\in [0,T)\times \R^d,$ the SDE
 \beq\label{Eab} \d X_{s,t}^x= b(t,X_{s,t}^x)\d t+ \si(t,X_{s,t}^x)\d W_t,\ \ \ t\in [s,T]\end{equation}
 is well-posed, so that  $(L_t^{a,b})_{t\in [s,T]}$ generates a unique diffusion process.   Moreover, for any $k\ge 1,$ there exists a constant $c(k)>0$ such that
 \beq\label{CK} \E \Big[\sup_{t\in [s,T]}|X_{s,t}^x|^k\Big] \le c(k)(1+|x|^k),\ \ (s,x)\in [0,T]\times\R^d .\end{equation}
 The associated Markov semigroup is given by
 $$P_{s,t}^{a,b}f(x):= \E[f(X_{s,t}^x)],\ \ \ \ 0\le s\leq t\le T, x\in\R^d, f\in \B_b(\R^d).$$ Since $(p_0,q_0)\in \scr K$, we have
 $$m_0:= \inf\Big\{m>1: \ff{(m-1)p_0}m \land \ff{(m-1)q_0}{ m} >1,\ \ff{dm}{p_0(m-1)}+\ff{2m}{q_0(m-1)}<2\Big\}\in (1,2).$$
 For a $\R^d\otimes\R^d$ valued differentiable function $a=(a^{ij})_{1\le i,j\le d}$, its divergence is an $\R^d$ valued  function defined as
$$\big({\rm div}a\big)^i:= \sum_{j=1}^d\pp_j a^{ij},\ \ \ 1\le i\le d.$$
   Our first result is the following.

 \beg{thm}\label{T1} Assume  $(A^{a,b})$ for $(a,b)=(a_i,b_i), i=1,2.$  Then for any $m\in(m_0,2)$, there exists a constant $c>0$ depending only on $m,K,d,T$  and $(p_i,q_i)_{0\le i\le l},$ such that
 for any $\psi\in \scr A$ and $\gg,\tt\gg\in \scr P,$
 \beq\label{PS1} \beg{split}& \W_\psi(P_{s,t}^{1,\gg},P_{s,t}^{2,\tt\gg})\le \ff{c\psi( (t-s)^{\ff 1 2})}{\ss{ t-s}} \W_1(\gg,\tt\gg)
  + c  \int_s^t \ff{\psi( (t-r)^{\ff 1 2}) \|a_1-a_2\|_{r,\infty}}{\ss{(r-s)(t-r) }}\d r \\
  &+c  \left(\int_s^t \left(\ff{\psi( (t-r)^{\ff 1 2}) \|a_1-a_2\|_{r,\infty}}{\ss{t-r} }\right)^{m}\d r\right)^{\frac{1}{m}}\\
  & +  c \int_s^t   \ff{\psi( (t-r)^{\ff 1 2}) } {\ss{t-r}} \big\{\|b_1-b_2\|_{r,\infty} +\|{\rm div}(a_1-a_2)\|_{r,\infty}\big\}\d r,\ \
     0\le s<t\le T, \  \gg,\tt\gg\in \scr P.\end{split} \end{equation}
 Moreover, for any $k\ge 1$,   there exists a
 constant $C>0$ depending only on $k,K,d,T$ and $(p_i,q_i)_{0\le i\le l},$ such that for any $\gg,\tt\gg\in \scr P$ and $0\le s\le t\le T,$
 \beq\label{PS2}  \W_k(P_{s,t}^{1,\gg},P_{s,t}^{2,\tt\gg})\le C\bigg[\W_k(\gg,\tt\gg) +   \int_s^t \|b_1-b_2\|_{r,\infty}\d r
   +   \bigg(\int_s^t \|a_1-a_2\|_{r,\infty}^{2}\d r\bigg)^{\ff 1 {2}}\bigg]. \end{equation}
\end{thm}

Next, we consider the following distribution dependent SDE on $\R^d$:
\beq\label{E1} \d X_t= b_t(X_t, \L_{X_t})\d t+  \si_t(X_t,\L_{X_t})\d W_t,\ \ t\in [0,T], \end{equation}
where   $\L_{X_t}$ is the distribution of $X_t$, and for some $k\ge 1$,
$$b: [0,T]\times\R^d\times\scr P_k \to\R^d,\ \ a: [0,T]\times\R^d\times\scr P_k  \to \R^d\otimes\R^d$$
are measurable, each  $a_t(x,\mu)$ is positive definite and $\si=\ss{2a}$.

Let $C^w_b([0,T];\scr P_k)$ be the set of
  all weakly continuous maps $\mu: [0,T]\to  \scr P_k$ such that
  $$\sup_{t\in [0,T]} \mu_t(|\cdot|^k)<\infty.$$
We call the SDE \eqref{E1}   well-posed for distributions in $\scr P_k$, if  for any initial value $X_0$ with $\L_{X_0}\in \scr P_k$ (correspondingly, any initial distribution $\nu\in \scr P_k$), the SDE has a unique solution (correspondingly,   a unique weak solution) with $(\L_{X_t})_{t\in [0,T]}\in C^w_b([0,T];\scr P_k).$
In this case, let
$P_t^*\nu:=\L_{X_t}$ for the solution with $\L_{X_0}=\nu$, and define
$$P_tf(\nu):= \int_{\R^d} f\d (P_t^*\nu),\ \ \nu\in \scr P_k, t\in [0,T], f\in \B_b(\R^d).$$
In particular, for $k=2$, the following log-Harnack inequality
\beq\label{LH} P_t\log f(\gg)\le \log P_tf(\tt\gg)+ \ff{c}t \W_2(\mu,\nu)^2,\ \ f\in \B_b^+(\R^d), t\in (0,T], \mu,\nu\in \scr P_2\end{equation} for some constant $c>0$ has been
established and applied  in  \cite{HW18, HW22a, RW, FYW1, FYW3}  for $\si_t(x,\mu)=\si_t(x)$ not dependent on $\mu$,   see also \cite{HRW19,HS,WZK} for extensions to the infinite-dimensional and    reflecting models.  When the noise coefficient is also distribution dependent and is $\W_2$-Lipschitz continuous, this inequality is established in the recent work \cite{23RW} by using a bi-coupling method.

In the following,  we consider more singular situation where $\si_t(x,\mu)$ may be not $\W_2$-Lipschitz continuous in $\mu$, and
the drift is singular in the time-spatial variables.
For any $\mu\in C_b^w([0,T];\scr P_k)$, let
$$a^\mu(t,x):= a_t(x,\mu_t),\ \ \ b^\mu(t,x):= b_t(x,\mu_t),\ \ \ t\in [0,T], x\in\R^d.$$
Correspondingly to $(A^{a,b})$, we make the following assumption.

\beg{enumerate} \item[$(B^{a,b})$] Let $k\in [1,\infty)$ and  $\psi\in \scr A$  with $\lim_{t\to0}\psi(t)=0$.
\item[$(1)$] $(A^{a,b})$ holds for $(a,b)= (a^\mu,b^\mu)$ uniformly in $\mu\in C_b^w([0,T];\scr P_k)$, with drift decomposition
$b^\mu=(b^\mu)^{(0)}+(b^\mu)^{(1)}$.
\item[$(2)$]    There exists    a constant $K>0$ such that
$$    \|a_t(\cdot,\gg)-a_t(\cdot,\tt\gg)\|_\infty
 \le K (\W_\psi+\W_k)(\gg,\tt\gg),\ \ \  t\in [0,T], \gg,\tt\gg\in \scr P_k.$$
\item[$(3)$]  There exist $p\ge 2$ and  $1\le \rr\in L^p([0,T])$, where $p=2$ if $\int_0^1 \ff{\psi(r)^2}r\d r<\infty$ and $p>2$ otherwise, such that for any $t\in [0,T]$ and $ \gg,\tt\gg\in \scr P_k,$
$$ \|b_t(\cdot,\gg)-b_t(\cdot,\tt\gg)\|_\infty+ \|{\rm div}( a_t(\cdot,\gg)-a_t(\cdot,\tt\gg))\|_\infty \le \rr_t (\W_\psi+\W_k)(\gg,\tt\gg).$$
\end{enumerate}

 \begin{rem}\label{EXa} We give a simple  example satisfying $(B^{a,b})$   for some  $\rr\in L^\infty([0,T])$, where $b$ contains a locally integrable term $b^{(0)}$,  and the dependence of $b$ and $\si$ in distribution is given by singular integral kernels.
 Let  $\psi\in \scr A$ with $\lim_{t\to0}\psi(t)=0$ and let
\beg{align*} &b_t(\cdot,\mu) = b_t^{(0)}+ \int_{\R^d} \tt b_t(\cdot, y) \mu(\d y),\\
&\si_t(\cdot,\mu)= \ss{\ll I+ \int_{\R^d} (\tt \si_t\tt\si_t^*)(\cdot,y) \mu(\d y)},\ \ (t,\mu)\in [0,T]\times \scr P_k,\end{align*}
where $\ll>0$ is a constant, $b^{(0)}: [0,T]\times\R^d\to\R^d$ satisfies $\|b^{(0)}\|_{\tt L_{q_0}^{p_0}(T)}<\infty$ for some $(p_0,q_0)\in \scr K$,
$\tt b: [0,T]\times\R^d\times\R^d\to\R^d$ is measurable such that
$$ |\tilde{b}_t(x,y)-\tilde{b}_t(\tilde{x},\tilde{y})|\leq K\big(|x-\tilde{x}|+ \psi(|y-\tilde{y}|)\big),\ \ x,\tilde{x},y,\tilde{y}\in\R^d,t\in[0,T]$$
holds for some constant $K>0$, and  $\tt \si: [0,T]\times\R^d\times\R^d\to \R^{d}\otimes\R^{d}$ is measurable and  bounded  such that
 \beg{align*}& \|\tilde{\sigma}_t(x,y)-\tilde{\sigma}_t(\tilde{x},\tilde{y})\|\leq K\big(|x-\tilde{x}| +\psi(|y-\tilde{y}|)\big),\\
 & |\nabla\tilde{\sigma}_t(\cdot,y)(x)-\nabla\tilde{\sigma}_t(\cdot,\tilde{y})(x)|\leq K \psi(|y-\tilde{y}|), \ \ x,\tilde{x},y,\tilde{y}\in\R^d,t\in[0,T].\end{align*}
 \end{rem}

We have the following result on the well-posedness and estimates on $(\W_\psi,\W_k)$  for $P_t^*$.

\beg{thm}\label{T2} Assume $(B^{a,b})$. Then the following assertions hold.
\beg{enumerate} \item[$(1)$]  The SDE $\eqref{E1}$ is well-posed for distributions in $\scr P_k$. Moreover, for any $n\in\mathbb N$,  there exists a constant $c >0$ such that any solution satisfies
\beq\label{ES1} \E\Big[\sup_{t\in [0,T]} |X_t|^n \Big|\F_0\Big] \le c (1+|X_0|^n).\end{equation}
\item[$(2)$] If  $\psi$ is a Dini function, i.e.
\beq\label{ESC} \int_0^1   \ff{ \psi(s)} s\,   \d s<\infty, \end{equation}
then   there exists a constant $c>0$ such that
\beq\label{ES4} \beg{split}&\W_\psi(P_t^*\gg,P_t^*\tt\gg)\le \ff{c\psi(t^{\ff 1 2})}{\ss t} \W_1(\gg,\tt\gg)+ c\W_k(\gg,\tt\gg),\\
&\W_k(P_t^*\gg,P_t^*\tt\gg)\le  c  \W_k(\gg,\tt\gg),\ \ \ \ t\in (0,T],\
\gg,\tt\gg\in \scr P_k.\end{split} \end{equation}\end{enumerate}
\end{thm}

\begin{rem}\label{EXa2} Theorem \ref{T2}(1) improves   existing well-posedness results for singular McKean-Vlasov SDEs  where  the coefficients are either   $(\W_k+\W_\aa)$-Lipschitz continuous in distribution  for some $\aa\in (0,1]$ and $ k\ge 1$ (see \cite{HWJMAA, HX23} and references therein), or satisfy some derivative conditions in distribution (see for instance \cite{CF}).\end{rem}

    To estimate $\W_{\psi}(P_t^*\gg,P_t^*\tt\gg)$ for worse $\psi $ not satisfying \eqref{ESC},  and to estimate the   relative entropy $\Ent(P_t^*\gg|P_t^*\tt\gg)$,   we need the drift to be Dini continuous in the spatial variable.

\beg{thm}\label{T3}  Assume $(B^{a,b})$ with $\|\rr\|_\infty<\infty$ and
  $\int_0^1\ff{\psi(r)^2}r\d r<\infty,$ and there exists  $\phi\in \scr A$ satisfying $\eqref{ESC}$ such that
  $$\sup_{\mu\in C_b^w([0,T];\scr P_k) } \big\{\|(b^\mu)^{(0)}\|_\infty+ [(b^\mu)^{0}]_\phi+\|\nn a^\mu \|_\infty\big\}<\infty.$$   Then the  following assertions hold.
 \beg{enumerate}\item[$(1)$ ] If $\psi(r)^2\log(1+r^{-1})\to 0$ as $r\to 0$, then there exists a constant $c>0 $ such that $\eqref{ES4}$ holds,  and for any
 $ t\in (0,T], \gg,\tt\gg\in \scr P_k,$
 \beq\label{LH1} \beg{split}& \Ent(P_t^*\gg|P_t^*\tt\gg)\le \ff{c\W_2(\gg,\tt\gg)^2}{  t}  \\
&\qquad +c\W_k(\gg,\tt\gg)^2\bigg(\ff 1 t\int_0^t \ff{\psi(r)^2}r\d r+ \ff{\psi(t^{\ff 1 2})^2}t  \log (1+t^{-1})\bigg). \end{split} \end{equation}
\item[$(2)$ ] If either $\|b\|_\infty<\infty$ or
\beq\label{E*}\sup_{(t,\mu)\in [0,T]\times\scr P_k} \big(\|\nn^i b_t(\cdot,\mu)\|_\infty + \|\nn^i\si_t(\cdot,\mu)\|_\infty\big)<\infty,\ \ i=1,2,\end{equation}
then there exists a constant $c>0$ such that $\eqref{ES4}$ holds,  and
\beq\label{LH2} \Ent(P_t^*\gg|P_t^*\tt\gg)\le \ff{c\W_2(\gg,\tt\gg)^2}{  t}
 +\ff{c\W_k(\gg,\tt\gg)^2} t\int_0^t \ff{\psi(r)^2}r\d r,\ \ t\in (0,T], \gg,\tt\gg\in \scr P_k.  \end{equation}\end{enumerate}
 \end{thm}
\begin{rem}\label{EXa3}  When $k\le 2$, $\eqref{LH} $  follows from  \eqref{LH2} or \eqref{LH1}.
  This improves \cite[Theorem 1.2]{23RW}, where  the $\W_2$-Lipschitz condition on the coefficients $(a,b)$ is relaxed as the  $(\W_\psi+\W_k)$-Lipschitz condition.\end{rem}

\


\section{Proof of Theorem \ref{T1}}

We first present a lemma to bound $\W_\psi$ by the total variation distance and $\W_1$.

\begin{lem}\label{LW} For any  $\psi\in \scr A, $
$$\W_{\psi}(\gamma,\tilde{\gamma})\leq \ss d\, \psi(\sqrt{t})\|\gamma-\tilde{\gamma}\|_{var}+\frac{d\psi(\sqrt{t})}{\sqrt{t}}\W_{1}(\gamma,\tilde{\gamma}),\ \ \gamma,\tilde{\gamma}\in\scr P_1. $$
\end{lem}

\beg{proof} Since $\psi$ is nonnegative and concave, we have
\beq\label{WF0} \psi(Rr)\le R\psi(r),\ \ \ r\ge 0, R\ge 1.\end{equation}
For any function $f$ on $\R^d$ with $[f]_\psi\le 1,$ let
$$f_t(x):= \E[f(x+B_t)],\ \ \ t\ge 0, x\in\R^d,$$
where $B_t$ is the standard Brownian motion  on $\R^d$ with $B_0=0$.   We have $\E[|B_t|^2]=dt.$  By $[f]_\psi\le 1$,   Jensen's inequality and \eqref{WF0},  we obtain
$$|f_t(x)-f(x)|\le \E[\psi(|B_t|)] \le \psi(\E|B_t|)\le \psi((dt)^{\ff 1 2})\le \ss d \psi(t^{\ff 1 2}),\ \ t\ge 0, x\in \R^d.$$
So,
\beq\label{1} \sup_{[f]_\psi\le 1} \big|\gg(f_t-f)-\tt\gg(f_t-f)\big|\le\ss d \psi(t^{\ff 1 2}) \|\gg-\tt\gg\|_{var},\ \ t\ge 0.\end{equation}
Next, for $[f]_\psi\le 1$, by Jensen's inequality, \eqref{WF0}, $\E|B_t|^2=dt$ and
$\E|B_t|\le \ss {dt},$
we obtain
\beg{align*} &|\nn f_t(x)|=\bigg|\nn_x \int_{\R^d} (2\pi t)^{-\ff d 2} \e^{-\ff{|x-y|^2}{2t}} (f(y)-f(z))\d y\bigg|_{z=x}\\
&\le (2\pi t)^{-\ff d 2} \int_{\R^d} \ff{|x-y|}t |f(y)-f(x)|\e^{-\ff{|x-y|^2}{2t}} \d y
\le \ff{1}t \E[|B_t|\psi(|B_t|)]\\
& \le \ff {\E|B_t|} t \psi\Big(\ff{\E|B_t|^2}{\E|B_t|}\Big)
 =\ff{\E|B_t|}{t} \psi\Big(\ff{(d\E|B_t|^2)^\ff 1 2}{\E|B_t|}  t^{\ff 1 2}\Big) \le  d  t^{-\ff 1 2}   \psi(t^{\ff 1 2}),\ \ t>0.\end{align*}
Combining this with \eqref{1} and noting that
$$\W_1(\gg,\tt\gg)=\sup_{\|\nn g\|\le 1} |\gg(g)-\tt\gg(g)|,$$ we derive that for any $f$ with $[f]_\psi\le 1$,
\beg{align*}&|\gg(f)-\tt\gg(f)|\le \big|\gg(f_t-f)-\tt\gg(f_t-f)\big|+ |\gg(f_t)-\tt\gg(f_t)|\\
&\le \ss d \psi(t^{\ff 1 2})\|\gg-\tt\gg\|_{var}+  dt^{-\ff 1 2} \psi(t^{\ff 1 2})\W_1(\gg,\tt\gg),\ \ t>0.\end{align*}
Then the proof is finished.
\end{proof}

Next, we present a gradient estimate on $P_{s,t}^{a,b}$.  All constants in the following   only depend  on $T,K,d$ and $(p_i,q_i)_{0\le i\le l}$.

\beg{lem} \label{LW2} Assume $(A^{a,b})$   without $\eqref{AA}$. Then there exists a constant $c>0$ such that for any $\psi\in \scr A,$
$$\sup_{[f]_\psi\le 1} \|\nn P_{s,t}^{a,b} f\|_\infty\le c(t-s)^{-\ff 1 2} \psi\big((t-s)^{\ff 1 2}\big),\ \ \ 0\le s<t\le T.$$
\end{lem}

\beg{proof}
(a) By \cite[Theorem 1.1]{XXZZ} or \cite[Theorem 2.1]{FYW3},  there exists a constant $c_1>0$ such that for any $0\le s<t\le T$ and $x\in\R^d$,  the Bismut formula
\beq\label{FF1} \nn P_{s,t}^{a,b}f(x)= \E\big[f(X_{s,t}^x) M_{s,t}^x\big]\end{equation}
holds for some random variable $M_{s,t}^x$ on $\R^d$ with
\beq\label{FF2} \E[M_{s,t}^x]=0,\ \ \ \E|M_{s,t}^x|^2\le c_1^2 (t-s)^{-1}. \end{equation}
So, for any $z\in \R^d$ and a function $f$ with $[f]_\psi\le 1$,
$$|\nn P_{s,t}^{a,b} f(x)| =\Big| \E\big[\{f(X_{s,t}^x) -f(z)\}M_{s,t}^x\big]\Big|\le \E\big[\psi(|X_{s,t}^x-z|)|M_{s,t}^x|\big].$$
By Jensen's inequality for the weighted probability $\ff {|M_{s,t}^x|\P}{\E |M_{s,t}^x|},$ we obtain
\beg{align*} &|\nn P_{s,t}^{a,b} f(x)|  \le \E[|M_{s,t}^x|] \psi\bigg(\ff{\E[| X_{s,t}^x-z|\cdot |M_{s,t}^x|]}{ \E[|M_{s,t}^x|]}\bigg) \\
&\le \E[|M_{s,t}^x|] \psi\bigg(\ff{ (\E[ |M_{s,t}^x|]^2)^{\ff 1 2}} {  \E[|M_{s,t}^x|]} \big(\E| X_{s,t}^x-z|^2\big)^{\frac{1}{2}}\bigg).\end{align*}
Combining this with \eqref{WF0} and \eqref{FF2}, we obtain
\beq\label{WF1}\sup_{[f]_\psi\le 1}  |\nn P_{s,t}^{a,b} f(x)|\le c_1 (t-s)^{-\ff 12}\inf_{z\in\R^d} \psi\Big( \big\{\E|X_{s,t}^x-z|^2\}^{\ff 1 2}\Big),\ \  0\le s<t\le T, x\in\R^d.\end{equation}

(b) To estimate $\inf_{z\in \R^d} \E|X_{s,t}^x-z|^2,$ we use Zvonkin's transform. By \cite[Theorem 2.1]{ZY}, there exist constants $\bb\in (0,1)$ and $\ll, C>0$ such that the PDE
\beq\label{PDE} (\pp_t +L_t^{a,b}-\ll)u_t= -b^{(0)}(t,\cdot),\ \ \ t\in [0,T], u_T=0\end{equation} for $u: [0,T]\times\R^d\to\R^d$ has a unique solution satisfying
\beq\label{WF2} \|u\|_\infty +\|\nn u\|_\infty+\sup_{x\ne y}\ff{|\nn u_t(x)-\nn u_t(y)|}{|x-y|^\bb} \le \ff 1 2,\end{equation}
\beq\label{WFM} \|\nn^2 u\|_{\tt L_{q_0}^{p_0}(T)} + \|(\pp_t+b^{(1)}\cdot\nn) u\|_{\tt L_{q_0}^{p_0}(T)}\le C.\end{equation}
By It\^o's formula, $Y_{s,t}:=\Theta_t(X_{s,t}^x), $ where $\Theta_t(y):= y+u_t(y),$  solves the SDE
$$\d Y_{s,t}= \bar b(t,Y_{s,t}) \d t +\bar\si(t,Y_{s,t})  \d W_t,\ \ t\in [s,T], Y_{s,s}= x+u_s(x),$$ where
\beq\label{BB} \bar b(t,\cdot):= (\ll u_t + b^{(1)})\circ \Theta_t^{-1},\ \ \bar\si(t,\cdot):=  \big\{(\nn \Theta_t)\si_t\big\}\circ \Theta_t^{-1}. \end{equation} By \eqref{WF2},
we find a constant $c_1>0$ such that
\beq\label{WF3} |\bar b(t,y)-\bar b(t,z)|\le c_1|y-z|,\ \ \ \|\bar\si(t,y)\|\le c_1,\ \ \ t\in [s,T], y,z\in\R^d.\end{equation}
Let
$$\ff{\d}{\d t} \theta_{s,t}=\bar b (t,\theta_{s,t})\big),\ \ \ t\in [s,T], \theta_{s,s} = Y_{s,s}= x+u_s(x).$$
By  It\^o's formula and \eqref{WF3}, we find a constant $c_2>0$ and a martingale $M_t$  such that
\beg{align*} \d |Y_{s,t}-\theta_{s,t}|^2 &= \Big\{2 \big\<Y_{s,t}-\theta_{s,t},  \bar b(t, Y_{s,t})  - \bar b(t,\theta_{s,t})\big\> + \|\bar \si (t,Y_{s,t})\|_{HS}^2 \Big\}\d t + \d M_t\\
&\le c_2 \big\{|Y_{s,t}-\theta_{s,t}|^2 +1\Big\}\d t+\d M_t,\ \ \ t\in [s,T], |Y_{s,s}-\theta_{s,s}|=0.\end{align*}
Thus,
$$\E\big[|Y_{s,t}-\theta_{s,t}|^2\big]\le c_2 \e^{c_2 T}(t-s),\ \ \ 0\le s\le t\le T.$$
Taking $z_{s,t}=\Theta_t^{-1}(\theta_{s,t})$ and noting that $\|\nn\Theta^{-1}\|_\infty<\infty$ due to $\|\nn u\|_\infty\le \ff 1 2$ in \eqref{WF2}, we find a constant $c_3>0$ such that
$$\E\big[|X_{s,t}^x-z_{s,t}|^2\big]= \E\big[|\Theta_t^{-1}(Y_{s,t})-\Theta_t^{-1}(\theta_{s,t})  |^2\big]\le c_3 (t-s),\ \ \ 0\le s\le t\le T.$$
Combining this with \eqref{WF1} and \eqref{WF0}, we finish the proof.
\end{proof}

Moreover, we estimate $\nn_y p_{s,t}^{a,b}(x,y)$, where $\nn_y$ is the gradient in $y$ and $p_{s,t}^{a,b}(x,\cdot)$ is the density function of $\L_{X_{s,t}^x}$. For any constant $\kk>0$, let
$$g_\kk(r,z):= (\pi\kk r)^{-\ff d 2} \e^{-\ff{|z|^2}{\kk r}},\ \ \ r>0, z\in\R^d$$ be the standard Gaussian heat kernel with parameter $\kk$.

\beg{lem} \label{LW3} Assume $(A^{a,b})$. Then   for any $m\in (m_0,2)$ there exists a constant $c(m)>0$ such that for any $t\in (0,T]$ and $0\le g_{\cdot,t}\in \B([0,t])$,
\beq\label{AVV}\begin{split} &\int_s^t \ff{g_{r,t} }{\ss {t-r}}  \d r   \int_{\R^d}  |\nn_y p_{s,r}^{a,b}(x,y) | \d y \\
&\le c(m)  \int_s^t \ff{g_{r,t}}{\ss{(t-r) (r-s)}} \d r+c(m) \bigg(\int_s^t\Big(\ff{g_{r,t} }{\ss {t-r}}\Big)^m\d r\bigg)^{\ff 1 m} ,
  \ \ s\in [0,t].
  \end{split}\end{equation} Consequently,
 there exists a constant $c >0$ such that
\beq\label{AVV'} \int_s^t  (t-r)^{-\ff 1 2}  \d r   \int_{\R^d}  |\nn_y p_{s,r}^{a,b}(x,y) | \d y \le c,
  \ \ 0\le s<t\le T.\end{equation}
\end{lem}

\beg{proof} Let $u_t$ be in \eqref{PDE}. By $(A^{a,b})$, $\si=\ss{2a}$,  \eqref{WF2} and \eqref{BB}, we  find a constant $c_1>0$ such that
$$|\bar b(t,x)-\bar b(t,y)|\le c_1 |x-y|,\ \ \|\bar \si (t,x)-\bar \si(t,y)\|\le c_1 |x-y|^{\aa\land \bb},\ \ \ t\in [0,T], x,y\in\R^d.$$
Let $\bar p_{s,t}(x,y)$ be the density function of $\L_{Y_{s,t}}$. According to \cite[Theorem 1.2]{MPZ}, there exists a constant $\kk\ge 1$ and some
$\theta_{s,t}:\R^d\to\R^d$ such that
\beq\label{WF4} |\nn_y^i \bar p_{s,t}(x,y)|\le \kk (t-s)^{-\ff i 2} g_\kk(t-s, \theta_{s,t}(x)-y),\ \ \ 0\le s<t\le T, x,y\in \R^d, i=0,1,\end{equation}
where $\nn^0 f:=f$.
Noting that $X_{s,t}^x=\Theta_t^{-1}(Y_{s,t})$, we have
 \begin{align}\label{PAP}p_{s,t}^{a,b}(x,y) = \bar p_{s,t} (\Theta_s(x),\Theta_t(y)) \big|{\rm det} (\nn\Theta_t(y))\big|.\end{align}
Combining this with \eqref{WF2}, \eqref{WF3} and \eqref{WF4}, we find a constant $c_2>0$ such that
 \beq\label{GRP}\begin{split} |\nn_y p_{s,t}^{a,b} (x,y)|\le &\,c_2  \kk (t-s)^{-\ff 1 2} g_\kk(t-s, \theta_{s,t}(\Theta_s(x))-\Theta_t(y))\big|{\rm det} (\nn\Theta_t(y))\big| \\
&+ c_2 \|\nn^2 u_t(y)\| p_{s,t}^{a,b}(x,y) ,\ \ 0\le s<t, x,y\in\R^d.
\end{split}\end{equation}
 Since $(p_0,q_0)\in \scr K,$ for any $m>m_0$, we have
\begin{align}\label{tpq}\tt p:= \ff{(m-1)p_0}m >1,\ \ \tt q:=  \ff{(m-1)q_0}m >1,\ \ \ff d{ \tt p}+\ff 2 {\tt q}<2.
\end{align}
By Krylov's estimate, see \cite[Theorem 3.1]{ZY}, we find a constant $c>0$ such that
\begin{equation}\label{kry}
\begin{split}&\int_s^t \d r  \int_{\R^d}  \|\nn^2 u_r(y)\|^{\ff m{m-1}} p_{s,r}^{a,b} (x,y) \d y\\
&= \E \int_s^t \|\nn^2 u_r\|^{\ff m{m-1}}(X_{s,r}^x)\d r\leq c \|\|\nn^2 u\|^{\ff m{m-1}}\|_{\tt L_{\tilde{q}}^{\tilde{p}}(s,t)}=c(\|\nn^2 u\|_{\tt L_{q_0}^{p_0}(s,t)})^{\frac{m}{m-1}}.
\end{split}\end{equation}
This together with \eqref{WFM}, \eqref{PAP} and \eqref{GRP} implies that for any $m\in(m_0,2)$, there exists a constant  $c(m) >0$ such that
 \beg{align*} & \int_s^t \ff{g_{r,t} }{\ss{t-r}} \d r \int_{\R^d}  |\nn_y p_{s,r}^{a,b} (x,y)| \d y
 \le c_2\kk \int_s^t g_{r,t}(t-r)^{-\ff 1 2}(r-s)^{-\ff 1 2} \d r \\
 &\qquad +
c_2  \bigg(\int_s^t \Big(\ff{g_{r,t} }{\ss{t-r}} \Big)^m  \d r\bigg)^{\ff 1 m} \bigg(\int_s^t \d r  \int_{\R^d}  \|\nn^2 u_r(y)\|^{\ff m{m-1}} p_{s,r}^{a,b} (x,y) \d y\bigg)^{\ff {m-1}m}
\\
& \le c(m)  \int_s^t \ff{g_{r,t}}{\ss{(t-r) (r-s)}} \d r  + c(m) \bigg(\int_s^t \Big(\ff{g_{r,t} }{\ss{t-r}} \Big)^m  \d r\bigg)^{\ff 1 m}.\end{align*}
So, \eqref{AVV} holds.  Letting $g_{r,t}\equiv 1$ and $m= \ff{m_0+2} 2$, we find a constant $c>0$ such that   \eqref{AVV} implies \eqref{AVV'}.

 \end{proof}

\beg{proof}[Proof of Theorem \ref{T1}]  By \eqref{MK}, it suffices to prove for   $\gg=\dd_x, \tt\gg=\dd_y, x,y\in\R^d$.

(a) We first consider $x=y$.
Let $f\in C_b^2(\R^d)$   with $[f]_\psi\le 1$. By It\^o's formula we have
$$P_{s,t}^{a_2,b_2}f(x)= f(x)+\int_s^t P_{s,r}^{a_2,b_2}(L_r^{a_2,b_2}f)(x)\d r,\ \ 0\le s\le t\le T.$$
This implies the Kolmogorov forward equation
\beq\label{KF} \pp_t P_{s,t}^{a_2,b_2}f= P_{s,t}^{a_2,b_2}(L_tf),\ \ {\rm a.e.}\ t\in [s,T].\end{equation}
On the other hand,   for $(p,q)\in \scr K$ and $t\in (0,T]$, let  ${\tt W}_{1,q,b_2^{(1)}}^{2,p}(0,t)$ be the set of all maps $u: [0,t]\times\R^d\to \R^d$ satisfying
$$\|u\|_{0\to t,\infty}+\|\nn u\|_{0\to t,\infty} +\|\nn^2 u\|_{\tt L_q^p(t)} +\|(\pp_s+b^{(1)}_2 \cdot\nn ) u\|_{\tt L_q^p(t)} <\infty.$$
By \cite[Theorem 2.1]{ZY},    the PDE
\beq\label{PDE1} (\pp_s+L_s^{a_2,b_2})u_s= - L_s^{a_2,b_2}f,\ \  s\in [0,t], u_t=0\end{equation}
has a unique solution in the class $ {\tt W}_{1,q,b_2^{(1)}}^{2,p}(0,t).$   So,  by  It\^o's formula   \cite[Lemma 3.3]{ZY},
 $$\d u_r(X_{s,r}^{2,x})= -L_r^{a_2,b_2} f(X_{s,r}^{2,x})+\d M_r,\ \ r\in [s,t]$$  holds for some martingale $M_r$.
This and \eqref{KF} yield
\beg{align*} & 0 = \E u_t(X_{s,t}^{2,x}) = u_s(x) - \int_s^t (P_{s,r}^{a_2,b_2} L_r^{a_2,b_2}f)\d r \\
&= u_s(x)- \int_s^t \ff{\d}{\d r} (P_{s,r}^{a_2,b_2}f)\d r= u_s(x)- P_{s,t}^{a_2,b_2}f + f,\ \ 0\le s\le t\le T.\end{align*}
Combining this with \eqref{PDE1},  we derive $P_{\cdot,t}^{a_2,b_2}f\in {\tt W}_{1,q,b_2^{(1)}}^{2,p}(0,t)$ for $t\in (0,T]$ and the Kolmogorov backward equation
 \beq\label{KB} \pp_sP_{s,t}^{a_2,b_2}f=\pp_s u_s = -L_s^{a_2,b_2}(u_s+f) =- L_s ^{a_2,b_2}P_{s,t}^{a_2,b_2}f,\ \ 0\le s\le t\le T.\end{equation}
  By It\^o's formula to  $P_{r,t}^{a_2,b_2} f(X_{s,r}^{1,x})$ for $r\in [s,t]$, see \cite[Lemma 3.3]{ZY}, we derive
\beg{align*}&P_{s,t}^{a_1,b_1}f(x)- P_{s,t}^{a_2,b_2} f(x)=\E \int_s^t \big(\pp_r+L_r^{a_1,b_1}\big) P_{r,t}^{a_2,b_2} f(X_{s,r}^{1,x})\d r \\
&= \int_s^t \d r\int_{\R^d} p_{s,r}^{a_1,b_1}(x,y)  \big(L_r^{a_1,b_1}-L_r^{a_2,b_2}\big) P_{r,t}^{a_2,b_2} f(y)  \d y.\end{align*}
By the integration by parts formula, we obtain
\beg{align*} &\bigg|\int_{\R^d} p_{s,r}^{a_1,b_1}(x,y) \big[ {\rm tr}\{(a_1-a_2)(r,y)\nn^2 P_{r,t}^{a_2,b_2} f (y)\}\big]\d y\bigg|\\
&= \bigg|\int_{\R^d}  \Big\<(a_1-a_2)(r,y) \nn_y p_{s,r}^{a_1,b_1}(x,y) + p_{s,r}^{a_1,b_1}(x,y){\rm div} (a_1-a_2)(r,y), \  \nn P_{r,t}^{a_2,b_2} f (y)\Big\>    \d y\bigg|.\end{align*}
Combining these with Lemma \ref{LW2}  and Lemma \ref{LW3},    for any $m\in (m_0,2)$, we find constants $c_1,c_2>0$ such that
\beg{align*} &|P_{s,t}^{a_1,b_1}f(x)- P_{s,t}^{a_2,b_2} f(x)|\le  c_1  \int_s^t \ff{\psi((t-r)^{\ff 1 2})\|a_1-a_2\|_{r,\infty}}{\ss{t-r}} \d r\int_{\R^d}  |\nn_y p_{s,r}^{a_1,b_1}(x,y)|\d y\\
&\qquad + c_1 \int_s^t \ff{\psi((t-r)^{\ff 1 2})}{(t-r)^{\ff 1 2}} \big(\|b_1-b_2\|_{r,\infty}+\|{\rm div}(a_1-a_2)\|_{r,\infty} \big)\d r \\
&\le c_2\int_s^t \ff{\psi((t-r)^{\ff 1 2})}{\ss{t-r}}\bigg( \ff{\|a_1-a_2\|_{r,\infty} }{\ss{r-s}}  + \|b_1-b_2\|_{r,\infty}+\|{\rm div}(a_1-a_2)\|_{r,\infty}  \bigg)\d r\\
& +c_2\left(\int_s^t \left(\ff{\psi((t-r)^{\ff 1 2})\|a_1-a_2\|_{r,\infty} }{\ss{t-r}}\right)^m\d r\right)^{\frac{1}{m}}=:I_{s,t}. \end{align*}
Therefore,
\beq\label{WT} \W_\psi\big(P_{s,t}^{1,x}, P_{s,t}^{2,x}\big)\le I_{s,t},\ \ \ 0\le s<t\le T, \ x\in\R^d.\end{equation}

(b) Let $x,y\in\R^d$ and $0\le s<t\le T.$  By the triangle inequality for $\W_\psi$, \eqref{WT}  and Lemma \ref{LW}, we obtain
\beq\label{*6}\beg{split}& \W_\psi(P_{s,t}^{1,x}, P_{s,t}^{2,y})\le \W_\psi(P_{s,t}^{1,x}, P_{s,t}^{2,x})+\W_\psi(P_{s,t}^{2,x}, P_{s,t}^{2,y})\\
&\le I_{s,t} + \psi\big((t-s)^{\ff 1 2}\big) \|P_{s,t}^{2,x}-P_{s,t}^{2,y}\|_{var} + \ff{\psi((t-s)^{\ff 1 2})}{\ss{t-s}} \W_1(P_{s,t}^{2,x}, P_{s,t}^{2,y}).\end{split}\end{equation}
By \cite[Theorem 2.1]{FYW3} or \cite[Theorem 1.1]{XXZZ}, $(A^{a,b})$ for $(a,b)=(a_2,b_2)$ implies that  for some   constant $c_3>0$,
$$\W_1(P_{s,t}^{2,x}, P_{s,t}^{2,y})\le c_3 |x-y|,\ \ \ \|P_{s,t}^{2,x}-P_{s,t}^{2,y}\|_{var}\le \ff{c_3 }{\ss{t-s}}|x-y|$$ holds for any $  0\le s<t\le T$ and $ x,y\in \R^d.$
Combining this with \eqref{*6}, we derive \eqref{PS1} for $\gg=\dd_x$ and $\tt\gg=\dd_y.$

(c)  It remains to prove \eqref{PS2}.  Let $u$ be in \eqref{PDE} for $(a,b)=(a_1,b_1)$. Let $\Theta_t(y):=y+u_t(y),$ and
$$Y_{s,t}^{1,x}= \Theta_t(X_{s,t}^{1,x}),\ \ \ Y_{s,t}^{2,y}= \Theta_t(X_{s,t}^{2,y}),\ \ \ t\in [s,T].$$
By It\^o's formula \cite[Lemma 3.3]{ZY}, we obtain
\beg{align*} &\d Y_{s,t}^{1,x} =  \big\{b_1^{(1)} (t,\cdot)+ \ll u_t\big\}(X_{s,t}^{1,x}) \d t + \big\{(\nn \Theta_t)\si_1(t,\cdot)\big\}(X_{s,t}^{1,x})\d W_t,\\
&\d Y_{s,t}^{2,y} =  \big\{b_1^{(1)} (t,\cdot)+ \ll u_t\big\}(X_{s,t}^{2,y}) \d t + \big\{(\nn \Theta_t)( b_2-b_1) + {\rm tr}[(a_2-a_1)(t,\cdot)\nn^2 u_t] \big\}(X_{s,t}^{2, y})\d t\\
&\qquad +   \big\{(\nn \Theta_t)\si_2(t,\cdot)\big\}(X_{s,t}^{2,y})\d W_t,\ \ \  t\in [s,T],\ Y_{s,s}^{1,x}=\Theta_s(x),\ Y_{s,s}^{2,y}=\Theta_s(y).\end{align*}
For any non-negative function $f$ on $\R^d$, let
$$\scr Mf(x):= \sup_{r\in (0,1]} \ff 1 {|B(x,r)|} \int_{B(x,r)} f(y)\d y,\ \ \ x\in \R^d, B(x,r):=\{y\in \R^d:|y-x|<r\}.$$
By $(A^{a,b})$  for $a=a_i$, $\si_i=\ss{2a_i},$ \eqref{WF2},  the maximal inequality in \cite[Lemma 2.1]{XXZZ},  and It\^o's formula, for any $k\ge 1$ we find a constant $c_1>1$ such that
\beq\label{90} c_1^{-1} |X_{s,t}^{1,x}-X_{s,t}^{2,y}|^{2k} \le  \xi_t:= |Y_{s,t}^{1,x}-Y_{s,t}^{2,y}|^{2k}\le c_1  |X_{s,t}^{1,x}-X_{s,t}^{2,y}|^{2k},\end{equation}
\beq\label{91} \d \xi_t\le c_1 \xi_t (1+\eta_t) \d t + c_1 \xi_t^{\ff {2k-1}{2k}}\gg_t \d t+ c_1\xi_t^{\ff{k-1}k} \|a_1-a_2\|_{t,\infty}^2\d t +\d M_t,\end{equation}
where $M_t$ is a  martingale  and
\beg{align*} &\gg_t:=\|b_1-b_2\|_{t,\infty} +\|a_1-a_2\|_{t,\infty}\|\nn^2 u_t\|(X_{s,t}^{2,y}),\\
&\eta_t:=       \scr M ( \|\nn\si_1 \|_{t,\infty}^2 +     \|\nn^2 u \|^2 ) (X_{s,t}^{1,x})+  \scr M(\|\nn\si_1 \|_{t,\infty}^2 +  \|\nn^2 u \|^2) (X_{s,t}^{2,y}).\end{align*}
 Note that for $q\in (\ff{2k-1}{2k},1)$,
\begin{align*}&\E\left\{\Big(\sup_{r\in[s,t]}\xi_r^q\Big)^{\ff {2k-1}{2kq}}\int_s^t\|a_1-a_2\|_{r,\infty}\|\nn^2 u_r\|(X_{s,r}^{2,y})\d r\right\}\\
&\leq \left(\E\sup_{r\in[s,t]}\xi_r^q\right)^{{\ff {2k-1}{2kq}}}\left(\E\left(\int_s^t\|a_1-a_2\|_{r,\infty}\|\nn^2 u_r\|(X_{s,r}^{2,y})\d r\right)^{\frac{2kq}{2kq-2k+1}}\right)^{\frac{2kq-2k+1}{2kq}}\\
&\leq \left(\E\sup_{r\in[s,t]}\xi_r^q\right)^{{\ff {2k-1}{2kq}}}\left(\int_s^t\|a_1-a_2\|^m_{r,\infty}\d r\right)^{\frac{1}{m}}\\
&\quad\quad\times\left(\E\left(\int_s^t\|\nn^2 u_r\|^{\frac{m}{m-1}}(X_{s,r}^{2,y})\d r\right)^{\frac{2(m-1)kq}{m(2kq-2k+1)}}\right)^{\frac{2kq-2k+1}{2kq}},\ \ m>1.
\end{align*}

 So, by the stochastic Grownwall inequality \cite[Lemma 2.8]{XZ} for $q\in (\ff{2k-1}{2k},1)$,
 \cite[Lemma 2.1]{XXZZ},  and the Krylov estimate  in \cite[Theorem  3.1]{ZY} which implies the Khasminskii inequality in \cite[Lemma 3.5]{XZ},
 we find  constants $c_2,c_3>0$ such that
\beg{align*}& \Big[\E\sup_{r\in [s,t]}\xi_r^q\Big]^{\ff 1 q} \le c_2 |x-y|^{2k} + c_2 \E\int_s^t \big\{\xi_r^{\ff {2k-1}{2k}}\gg_r\d r+\xi_r^{\ff{k-1}k}\|a_1-a_2\|_{r,\infty}^2\big\}\d r\\
 &\le c_2|x-y|^{2k} + c_2 \E\bigg[\Big(\sup_{r\in [s,t]}\xi_r^q\Big)^{\ff{2k-1}{2kq}} \int_s^t \gg_r\d r + \Big(\sup_{r\in [s,t]}\xi_r^q\Big)^{\ff{k-1}{kq}} \int_s^t \|a_1-a_2\|_{r,\infty}^2\d r\bigg]\\
 &\le c_2 |x-y|^{2k}+\ff 1 2 \Big[\E\sup_{r\in [s,t]}\xi_r^q\Big]^{\ff 1 q} + c_3 \bigg(\int_s^t \|a_1-a_2\|_{r,\infty}^2\d r\bigg)^k +
   c_3 \bigg(\int_s^t \|b_1-b_2\|_{r,\infty}\d r\bigg)^{2k}  \\
&+ c_3   \bigg(\int_s^t \|a_1-a_2\|_{r,\infty}^m \d r\bigg)^{\ff{2k}m}    \left(\E\left(\int_s^t\|\nn^2 u_r\|^{\frac{m}{m-1}}(X_{s,r}^{2,y})\d r\right)^{\frac{2(m-1)kq}{m(2kq-2k+1)}}\right)^{\frac{2kq-2k+1}{q}},\ \ m>1.
 \end{align*}  Noting that \cite[Theorem 2.1(3)]{Ren} implies
 $$\Big[\E\sup_{r\in [s,t]}\xi_r^q\Big]<\infty,$$
 we obtain
 \beq\label{92} \beg{split} &\Big[\E\sup_{r\in [s,t]}\xi_r^q\Big]^{\ff 1 q} \le 2c_2 |x-y|^{2k}  + 2c_3 \bigg(\int_s^t \|a_1-a_2\|_{r,\infty}^2\d r\bigg)^k \\
 &\quad +
   2c_3 \bigg(\int_s^t \|b_1-b_2\|_{r,\infty}\d r\bigg)^{2k}   \\
   &\quad   + 2c_3   \bigg(\int_s^t \|a_1-a_2\|_{r,\infty}^m \d r\bigg)^{\ff{2k}m}     \left(\E\left(\int_s^t\|\nn^2 u_r\|^{\frac{m}{m-1}}(X_{s,r}^{2,y})\d r\right)^{\frac{2(m-1)kq}{m(2kq-2k+1)}}\right)^{\frac{2kq-2k+1}{q}}.
 \end{split}  \end{equation}
Recall that $(\tilde{p},\tilde{q})$ is defined in \eqref{tpq}. By  \eqref{WFM},   \cite[Theorem  3.1]{ZY}  and  \cite[Lemma 3.5]{XZ} , we find a constant $c_4>0$ such that
 \begin{align*}  &\E\left(\int_s^t\|\nn^2 u_r\|^{\frac{m}{m-1}}(X_{s,r}^{2,y})\d r\right)^{\frac{2(m-1)kq}{m(2kq-2k+1)}}\\
&\le c_4 (\|\|\nn^2 u\|^{\ff m{m-1}}\|_{\tt L_{\tt q}^{\tt p}(s,t)})^{\frac{2(m-1)kq}{m(2kq-2k+1)}} = c_4(\|\nn^2 u\|_{\tt L_{q_0}^{p_0}(0,T)})^{\frac{2kq}{2kq-2k+1}} <\infty.
\end{align*}
 Combining this with   \eqref{92},  we find a constant $c_5>0$ such that
\beg{align*} &\big(\E|Y_{s,t}^{1,x}-Y_{s,t}^{2,y}|^k\big)^2 \le  \Big[\E\sup_{r\in [s,t]}\xi_r^q\Big]^{\ff 1 q} \le c_5 |x-y|^{2k} + c_5 \bigg(\int_s^t \|b_1-b_2\|_{r,\infty}\d r\bigg)^{2k} \\
&\qquad + c_5 \bigg(\int_s^t \|a_1-a_2\|_{r,\infty}^m\d r\bigg)^{\ff{2k}m} +c_5 \bigg(\int_s^t \|a_1-a_2\|_{r,\infty}^2\d r\bigg)^{k}.\end{align*}
Noting that   \eqref{90} implies
$$\W_k(P_{s,t}^{1,x}, P_{s,t}^{2,y})^k\le \ss{c_1} \E|Y_{s,t}^{1,x}-Y_{s,t}^{2,y}|^k,$$
by Jensen's inequality we derive  \eqref{PS2} for some constant $C>0$ and $\gg=\dd_x,\tt\gg=\dd_y.  $

\end{proof}

\section{Proof of Theorem \ref{T2}}

Once the well-posedness of \eqref{E1} is proved,     the proof of \cite[(1.5)]{HWJMAA} implies   \eqref{ES1} under $(B^{a,b})$. We skip the details to save space.
So, in the following we only prove the well-posedness and  estimate \eqref{ES4}.

(a) Let $X_0$ be $\F_0$-measurable with $\gg:=\L_{X_0}\in \scr P_k.$ Let
$$\C_{T}^\gg:= \big\{\mu\in C([0,T]; \scr P_k):\ \mu_0=\gg\big\}.$$
For any $\ll\ge 0,$ $C_{T}^\gg$ is a complete space under the metric
$$\rr_\ll(\mu, \tt\mu):=\sup_{t\in [0,T]} \e^{-\ll t} \big\{\W_\psi(\mu_t,\tt\mu_t)+\W_k(\mu_t,\tt\mu_t)\big\}.$$

For any   $\mu\in C([0,T];\scr P_k)$, let
$$b_t^{\mu}(x):=b_t(x,\mu_t),\ \ \sigma_t^{\mu}(x)=\sigma_t(x,\mu_t),\ \ \ (t,x)\in [0,T]\times \R^d.$$
According to \cite[Theorem 2.1]{Ren},   $(B^{a,b})$ implies that the SDE
$$\d X_t^\mu= b_t^\mu(X_t^\mu)\d t+ \si_t^\mu(X_t^\mu)\d W_t,\ \ \ t\in [0,T], X_0^\mu=X_0$$
is well-posed, and
$$\E\Big[\sup_{s\in [0,T]} |X_t^\mu|^k \Big]<\infty. $$
So,  we define a map
$$\Phi^\gg: \C_T^\gg\to \C_T^\gg;\ \
 \mu \mapsto \big\{(\Phi^\gg \mu)_t:=\L_{X_t^\mu}\big\}_{t\in [0,T]}.$$
According to \cite[Theorem 3.1]{HRW}, if $\Phi^\gg$ has a unique fixed point in $\C_T^\gg$, then \eqref{E1} is well-posed for distributions in $\scr P_k.$

(b) Let $\tt\gg\in \scr P_k$ which may be different from $\gg$, and let $\tt \mu\in \C_T^{\tt\gg}$.  We estimate the $\rr_\ll$-distance between $\Phi^\gg\mu$ and $
\Phi^{\tt\gg}\tt\mu.$
  By Theorem \ref{T1} and $(B^{a,b})$, for any $m \in(m_0,2)$, there exist constants $c_1,c_2>0$ such that
  \beg{align*} & \W_\psi\big((\Phi^\gg\mu)_t, (\Phi^{\tt\gg}\tt\mu)_t\big)+  \W_k\big((\Phi^\gg\mu)_t, (\Phi^{\tt\gg}\tt\mu)_t\big)\\
&\le \ff{c_1 \psi(t^{\ff 1 2})}{\ss t} \W_k(\gg,\tt\gg) +   c_1 \bigg(\int_0^t \|a^\mu-a^{\tt\mu}\|_{r,\infty}^2\d r\bigg)^{\ff 1 2}\\
&\quad+c_1  \left(\int_0^t \left(\ff{\psi( (t-r)^{\ff 1 2}) \|a^\mu-a^{\tilde{\mu}}\|_{r,\infty}}{\ss{t-r} }\right)^{ m}\d r\right)^{\frac{1}{ m}}\\
&\qquad + c_1 \int_0^t  \ff{c_1 \psi((t-r)^{\ff 1 2})}{\ss {t-r}}  \Big(\ff{\|a^\mu-a^{\tt\mu}\|_{r,\infty} }{\ss r} +\|b^\mu-b^{\tt\mu}\|_{r,\infty}
+\|{\rm div}(a^\mu-a^{\tt\mu})\|_{r,\infty}\Big)\d r \\
&\le \ff{c_1 \psi(t^{\ff 1 2})}{\ss t} \W_k(\gg,\tt\gg) + c_2 \bigg(\int_0^t  \Big(\W_\psi(\mu_r,\tt\mu_r)+\W_k(\mu_r,\tt\mu_r)\Big)^2 \d r \bigg)^{\ff 1 2}\\
&\quad+c_2  \left(\int_0^t \left(\ff{\psi( (t-r)^{\ff 1 2}) (\W_\psi(\mu_r,\tt\mu_r)+\W_k(\mu_r,\tt\mu_r))}{\ss{t-r} }\right)^{ m}\d r\right)^{\frac{1}{m}}\\
&\qquad + c_2 \int_0^t \ff{ \psi((t-r)^{\ff 1 2})}{\ss {r(t-r)}} (1+\sqrt{r}\rho_r)\Big(\W_\psi(\mu_r,\tt\mu_r)+\W_k(\mu_r,\tt\mu_r)\Big) \d r \\
\end{align*}
Let $\gg=\tt\gg$. We obtain
$$\rr_\ll(\Phi^\gg\mu,\Phi^{\gg}\tt\mu) \le \dd(\ll) \rr_\ll(\mu,\tt\mu),$$
where by $(B^{a,b})$ and $m\in  (m_0, 2)$,  as $\ll\to\infty$ we have
\begin{align*}\dd(\ll)&:= c_2\sup_{t\in [0,T]} \bigg[\int_0^t \ff{ \psi((t-r)^{\ff 1 2})\e^{-\ll(t-r)}}{\ss {t-r}}\Big(\ff 1 {\ss r}+\rr_r\Big)\d r + \bigg(\int_0^t\e^{-2\ll(t-r)}\d r\bigg)^{\ff 1 2} \bigg]\\
& +c_2  \left(\int_0^t \left(\ff{\psi( (t-r)^{\ff 1 2}) \e^{-\lambda(t-r)}}{\ss{t-r} }\right)^{ m}\d r\right)^{\frac{1}{ m}}\to 0.
\end{align*}
So, $\Phi^\gg$ is $\rr_\ll$-contractive on $\C_T^\gg$ for large $\ll>0$, and hence  has a unique fixed point. This implies the well-posedness of \eqref{E1}
for distributions in $\scr P_k$.

(c)  For $s\in [0,T),$  let $P_{s,t}^*\gg=\L_{X_{s,t}^\gg}, $ where $X_{s,t}^\gg$ solves \eqref{E1} for $t\in [s,T]$ and $\L_{X_{s,s}^\gg}=\gg.$
 By \eqref{ES1} for $s$ replacing $0$, we have
$$\sup_{t\in [s,T]} (P_{s,t}^*\gg)(|\cdot|^k)<\infty,\ \ \ \gg\in \scr P_k.$$ Since $\psi$ has growth slower than linear, and  \eqref{WF0} implies the boundedness of $\ff r {\psi(r)}$  for $r\in [0,T],$ this implies that  for any $\gg,\tt\gg\in \scr P_k$ and $s\in [0,T),$
\beq\label{GGT}  \sup_{r\in [s,t] }  (\W_\psi+\W_k)(P_{s,r}^*\gg,P_{s,r}^*\tt\gg) <\infty,\ \ t\in [s,T],
\end{equation}
\beq\label{GGT'}\GG_{s,t}:=  \sup_{r\in [s,t] }  \ff{\ss{r-s}}{\psi((r-s)^{\ff 1 2})} (\W_\psi+\W_k)(P_{s,r}^*\gg,P_{s,r}^*\tt\gg)<\infty,\ \ t\in [s,T].\end{equation}
Let
\beg{align*} &a_1(t,x):= a_t(x,P_{s,t}^*\gg),\ \ b_1(t,x):= b_t(x,P_{s,t}^* \gg),\\
&a_2(t,x):= a_t(x,P_{s,t}^*\tt\gg),\ \ b_1(t,x):= b_t(x,P_{s,t}^*\tt\gg),\ \ (t,x)\in [s,T]\times \R^d.\end{align*}
Then  $P_{s,t}^*\gg=P_{s,t}^{1,\gg}, P_{s,t}^*\tt\gg=P_{s,t}^{2,\tt\gg},$ and \eqref{MK}  implies
\beq\label{MK'} P_{s,t}^*\gg=\int_{\R^d} P_{s,t}^{1,x}\gg(\d x),\ \ \
P_{s,t}^*\tt\gg=\int_{\R^d} P_{s,t}^{2,x}\tt\gg(\d x).\end{equation}
Thus,  by  Theorem \ref{T1} and $(B^{a,b})$, for  any $ m \in(m_0,2)$, we find a constant $k_0>0$ such that
\beq\label{A1} \beg{split} &\W_\psi(P_{s,t}^*\gg, P_{s,t}^*\tt\gg)= \W_\psi(P_{s,t}^{1,\gg}, P_{s,t}^{2,\tt\gg})
 \le \ff{k_0 \psi((t-s)^{\frac{1}{2}})}{\ss{t-s}}\W_1(\gg,\tt\gg)\\
 &+ k_0\int_s^t \ff{ \psi((t-r)^{\ff 1 2})}{ \ss{(r-s)(t-r)}} \Big(1+\rr_{r}\ss{r-s}\Big)\big(\W_\psi+\W_k\big)(P_{s,r}^*\gg,P_{s,r}^*\tt\gg)\d r\\
 &+k_0  \left(\int_s^t \left(\ff{\psi( (t-r)^{\ff 1 2}) \big(\W_\psi+\W_k\big)(P_{s,r}^*\gg,P_{s,r}^*\tt\gg)}{\ss{t-r} }\right)^{ m}\d r\right)^{\frac{1}{ m}},\end{split}\end{equation}
\beq\label{A2} \beg{split} &\W_k(P_{s,t}^*\gg, P_{s,t}^*\tt\gg)= \W_k(P_{s,t}^{1,\gg}, P_{s,t}^{2,\tt\gg})
 \le k_0  \W_k(\gg,\tt\gg) \\
 & +k_0\int_s^t \rr_r  \big(\W_\psi+\W_k\big)(P_{s,r}^*\gg,P_{s,r}^*\tt\gg)\d r+k_0 \bigg( \int_s^t \big(\W_\psi+\W_k\big)^2(P_{s,r}^*\gg,P_{s,r}^*\tt\gg)\d r \bigg)^{\ff 1 2}.\end{split}\end{equation}
By combining these with the definition of $\GG_{s,t}$ in \eqref{GGT'}, we find a constant $k_1>0$ such that
\beq\label{GGT2} \beg{split}&\GG_{s,t} \le k_1  \W_k(\gg,\tt\gg)  + k_1 \GG_{s,t} h(t-s),\ \ 0\le s<t\le T,\\
&h(t):=  \sup_{(s,\theta)\in (0,t]\times [0,T-t]}   \ff{\ss s}{\psi(s^{\ff 1 2})} \int_0^s \ff{\psi(r^{\ff 1 2})\psi((s-r)^{\ff 1 2})}{ \ss{r(s-r)}} \Big(\ff 1 {\ss r}+\rr_{\theta+r}\Big)\d r \\
&\qquad +\sup_{s\in (0,t]}   \ff{\ss s}{\psi(s^{\ff 1 2})}  \left(\int_0^s \left(\ff{\psi( (s-r)^{\ff 1 2})\psi(r^{\frac{1}{2}})}{\ss r\ss{s-r} }\right)^{ m}\d r\right)^{\frac{1}{ m}}\\
&\qquad +  \bigg( \int_0^t \Big(\ff{\psi(r^{\ff 1 2})}{\ss r}\Big)^2\d r\bigg)^{\ff 1 2},\ \ t\in (0,T].\end{split} \end{equation}
Note that
\beg{equation}\begin{split} \label{BB1}&\ff{\ss s}{\psi(s^{\ff 1 2})} \int_0^s \ff{\psi(r^{\ff 1 2})\psi((s-r)^{\ff 1 2})}{r\ss{s-r}} \d r\\
&\leq\ff{\ss s}{\psi(s^{\ff 1 2})}\bigg(\int_0^{\ff s 2} \ff {\psi(s^{\ff 1 2})}{\ss{s/2}}\cdot\ff{\psi(r^{\ff 1 2})}{r}\d r +\int_{\ff s 2}^s \ff{\psi((s-r)^{\ff 12})}{s-r}\cdot\ff{\ss s \psi(s^{\ff 1 2})}{s/2}\d r \bigg)\\
&\le \big(2+\ss 2 \big)\int_0^{\ff s 2} \ff{\psi(r^{\ff 1 2})}r\d r=  2\big(2+\ss 2\big)  \int_0^{\ss{s/2}}\ff{\psi(r)}r\d r.
\end{split}\end{equation}
Similarly, we have
 {\beg{equation}\begin{split} \label{BB2}& \ff{\ss s}{\psi(s^{\ff 1 2})}  \left(\int_0^s \left(\ff{\psi( (s-r)^{\ff 1 2})\psi(r^{\frac{1}{2}})}{\ss r\ss{s-r} }\right)^{ m}\d r\right)^{\frac{1}{ m}} \\
&\leq\sqrt{2}\left(\left(\int_0^{\ff s 2}\left( \ff{\psi(r^{\ff 1 2})}{\ss r}\right)^{ m}\d r\right)^{^{\frac{1}{  m}}} +\left(\int_{\ff s 2}^s \left(\ff{\psi((s-r)^{\ff 12})}{s-r}\right)^{ m}\d r \right)^{\frac{1}{ m}}\right)\\
&\le 2\sqrt{2}\left(\int_0^{\ff s 2}\left( \ff{\psi(r^{\ff 1 2})}{\ss r}\right)^{ m}\d r\right)^{^{\frac{1}{ m}}},\end{split}\end{equation}
  \beg{equation}\begin{split} \label{BB3}&\ff{\ss s}{\psi(s^{\ff 1 2})} \int_0^s \ff{\psi(r^{\ff 1 2})\psi((s-r)^{\ff 1 2})}{\ss{r(s-r)}}\rr_{\theta+r} \d r\\
&= \ff{\ss s}{\psi(s^{\ff 1 2})}\bigg(\int_0^{\ff s 2} \ff {\psi(s^{\ff 1 2})}{\ss{s/2}}\cdot\ff{\psi(r^{\ff 1 2})}{\ss r}\rr_{\theta+r}\d r
+\int_{\ff s 2}^s \ff{\psi((s-r)^{\ff 12})}{\ss{s-r}}\cdot\ff{\ss s\psi(s^{\ff 1 2})}{s/\ss 2}\rr_{\theta+r}\d r \bigg)\\
&\le 2\ss 2 \int_0^s \Big(\ff{\psi(r^{\ff 1 2})}{\ss r} + \ff{\psi((s-r)^{\ff 1 2})}{\ss{s-r}}\Big)\rr_{\theta+r} \d r
\le  4 \ss 2\bigg( \int_0^s \ff{\psi(r^{\ff 1 2})^2}{  r} \d r\bigg)^{\ff 1 2}
\bigg(\int_0^T \rr_r^2\d r\bigg)^{\ff 1 2}.\end{split}\end{equation}
Combining these with    \eqref{ESC}, we  conclude that $h(t)$ defined in  \eqref{GGT2} satisfies $h(t)\to 0$ as $t\to 0$. Letting $r_0>0$ such that $k_1h(t)\le \ff 1 2$ for $t\in [0,r_0],$ we deduce form \eqref{GGT'} and  \eqref{GGT2} that
\beg{align*} \ff{\ss{t-s}}{\psi((t-s)^{\ff 1 2})} (\W_\psi+\W_k)(P_{s,t}^*\gg,P_{s,t}^*\tt\gg)
\le \GG_{s,t}
  \le 2 k_1 \W_k(\tt\gg,\gg)\end{align*} holds for all $  s\in [0,T)$ and $t\in (s, (s+r_0)\land T].$  Consequently,
 \beg{align*}&(\W_\psi+\W_k)(P_{s,t}^*\gg,P_{s,t}^*\tt\gg)\le \ff{2k_1\psi((t-s)^{\ff 1 2})} {\ss{t-s}}\W_k(\gg,\tt \gg),\\
 &\qquad  \ s\in [0,T), t\in (s, (s+r_0)\land T],\ \gg,\tt\gg\in \scr P_k.\end{align*}
 Combining this with the flow property
$$P_{s,t}^*= P_{r,t}^*P_{s,r}^*,\ \ \ 0\le s\le r\le t\le T,$$
we find a constant $k_2>0$ such that
\beq\label{G*} (\W_\psi+\W_k)(P_{s,t}^*\gg, P_{s,t}^*\tt\gg)
\le \ff{k_2\psi((t-s)^{\ff 1 2})}{\ss{t-s}} \W_k(\gg,\tt\gg),\ \ t\in (s,T], \gg,\tt\gg\in \scr P_k.\end{equation}
By the conditions on $\psi$ in $(B^{a,b})(3)$ and \eqref{ESC}, we have
 \begin{align*}&\sup_{t\in (0,T]}\Bigg\{ \int_0^t\ff{\psi(r^{\ff 1 2})\psi((t-r)^{\ff 1 2})}{r\ss{t-r}}\Big(1+\rr_r\ss r\Big)\d r+ \left(\int_0^t\Big(\ff{\psi(r^{\ff 1 2})}{\ss r}\Big)^2\d r\right)^{\frac{1}{2}}\\
&\qquad\quad+\left(\int_0^t \left(\ff{\psi( (t-r)^{\ff 1 2})\psi(r^{\frac{1}{2}})}{\ss r\ss{t-r} }\right)^{m}\d r\right)^{\frac{1}{m}}\Bigg\}<\infty.
\end{align*}

Therefore, substituting  \eqref{G*}  into \eqref{A1} and \eqref{A2}, we derive   \eqref{ES4} for some constant $c>0$.

\section{Proof of Theorem \ref{T3}}

(a) We use the notations in step (c) in the proof of Theorem \ref{T2}.
By Pinsker's inequality, \cite[(1.3)]{23RW} and  $(B^{a,b})$ with $\|\rr\|_\infty<\infty$, we find   constants $\vv\in (0,\ff 1 2], c_1>0$ such that
\beg{align*} &\|P_{s,t}^{1,x}-P_{s,t}^{2,y}\|_{var} \le \ss{2 \Ent(P_{s,t}^{1,x}|P_{s,t}^{2,y})} \\
&\le \ff{c_1|x-y|}{\ss{t-s}} +\ff{c_1}{\ss{t-s}}\bigg(\int_s^t(\W_\psi+\W_k)^2(P_{s,r}^*\gg, P_{s,r}^*\tt\gg)\d r\bigg)^{\ff 1 2} \\
&\quad + c_1 \ss{\log(1+(t-s)^{-1})} \sup_{r\in [s+\vv(t-s),t]}
(\W_\psi+\W_k)^2(P_{s,r}^*\gg, P_{s,r}^*\tt\gg)\d r\bigg),\ \ t\in [s,T].\end{align*}
Combining this with \eqref{MK'} and Lemma \ref{LW}, we obtain
\beq\label{GP1} \beg{split} & \W_\psi(P_{s,t}^*\gg, P_{s,t}^*\tt\gg) - \ff{\psi((t-s)^{\ff 1 2})}{\ss{t-s}} \W_1(P_{s,t}^{1,\gg}, P_{s,t}^{2,\tt\gg})\le   \psi((t-s)^{\ff 1 2}) \|P_{s,t}^{1,\gg}-P_{s,t}^{2,\tt\gg}\|_{var}\\
&\le \ff{\psi((t-s)^{\ff 1 2})}{\ss{t-s}} \bigg(\int_s^t(\W_\psi+\W_k)^2(P_{s,r}^*\gg, P_{s,r}^*\tt\gg)\d r\bigg)^{\ff 1 2}\\
&\quad + c_1 \psi((t-s)^{\ff 1 2}) \ss{\log(1+(t-s)^{-1})} \sup_{r\in [s+\vv(t-s),t]}
(\W_\psi+\W_k) (P_{s,r}^*\gg, P_{s,r}^*\tt\gg) \end{split}
\end{equation} for $ t\in [s,T].$
On the other hand, since $b^{(0)}$ is bounded, $\|b^{(0)}\|_{\tt L_{q_0}^{p_0}(T)}<\infty$ holds for any $p_0,q_0>2$, so that \eqref{PS2} holds for $m=2$. Then there exists a constant $c_2>0$ such that
\beq\label{GP2}\beg{split}& \W_1(P_{s,t}^{1,\gg}, P_{s,t}^{2,\tt\gg})\le \W_k(P_{s,t}^{1,\gg}, P_{s,t}^{2,\tt\gg})\\
&\le c_2 \W_k(\gg,\tt\gg) +c_2 \bigg(\int_s^t(\W_\psi+\W_k)^2(P_{s,r}^*\gg, P_{s,r}^*\tt\gg)\d r\bigg)^{\ff 1 2}.\end{split}\end{equation}
Combining this with \eqref{GP1},   we find a constant $c_3>0$ such that instead of \eqref{GGT2} we have
\beq\label{LST} \beg{split} & \GG_{s,t}\le c_3 \W_k(\gg,\tt\gg) + c_2 h(t-s) \GG_{s,t},\ \ \ 0\le s\le t\le T,\\
&h(t):= \bigg(\int_0^t \ff{\psi(s^{\ff 1 2})^2}s\d s\bigg)^{\ff 1 2}+ \sup_{r\in (0,t]} \psi(r^{\ff 12})\ss{\log(1+r^{-1})},\ \ t>0.\end{split}\end{equation}
Since  $\int_0^1\ff{\psi(r)^2}r\d r<\infty$,  we have
$h(t)\to 0$ as $t\to 0$ if $\lim_{r\to 0} \psi(r)^2\log(1+r^{-1})=0,$   so that  \eqref{ES4} follows as explained in step (c) in the proof of Theorem \ref{T2}.

(b) Next, by \eqref{MK'}, \cite[(1.3)]{23RW} and  $(B^{a,b})$ with $\|\rr\|_\infty<\infty, $ we find   constants $\vv\in (0,\ff 1 2], c_1>0$ such that for any $\gg,\tt\gg\in \scr P_k,$
\beg{align*} &\Ent(P_t^*\gg|P_t^*\tt\gg)\le \ff{\W_2(\gg,\tt\gg)^2}t +\ff{c_1}t \int_0^t (\W_\psi+\W_k)^2(P_{r}^*\gg, P_{r}^*\tt\gg)\d r\\
&\qquad + c_1 \log(1+t^{-1}) \sup_{r\in [\vv t,t]} (\W_\psi+\W_k)^2(P_{r}^*\gg, P_{r}^*\tt\gg),\ \ t\in (0,T].\end{align*}
Combining this with \eqref{ES4}, we find a constant $c>0$ such that \eqref{LH1} holds.

(c) If either $\|b\|_\infty<\infty$ or \eqref{E*} holds, then we may apply \cite[(1.4)]{23RW}
to delete the term $\log(1+(t-s)^{-1})$ from the above calculations, so that
$h(t)$ in \eqref{LST}  becomes $ \big(\int_0^t \ff{\psi(s^{\ff 1 2})^2}s\d s\big)^{\ff 1 2}$ which goes to $0$ as $t\to 0$. Therefore,
\eqref{ES4} and \eqref{LH2} hold for some constant $c>0$ as shown above.

\end{document}